\documentclass[letterpaper, 10 pt, conference]{ieeeconf}  
\usepackage[letterpaper, top=57pt, left=48pt, right=54pt, bottom=45pt]{geometry}

\IEEEoverridecommandlockouts                              

\usepackage{graphicx} 
\usepackage{amsmath} 
\usepackage{amssymb}  
\usepackage{dsfont}
\usepackage[linesnumbered,ruled]{algorithm2e}
\usepackage{booktabs}
\usepackage{scalerel}
\usepackage[font=scriptsize]{caption}
\usepackage{tikz}
\usepackage{xcolor} 
\usepackage{subcaption}
\usepackage{url}
\usepackage[pdfstartview=FitV, urlcolor=blue]{hyperref}

\newcommand{\vx}{\mathbf{x}}
\newcommand{\vy}{\mathbf{y}}
\newcommand{\va}{\mathbf{a}}

\newcommand{\A}{\mathcal{A}}

\newcommand{\D}{\mathcal{D}}
\newcommand{\E}{\mathbb{E}}
\newcommand{\R}{\mathbb{R}}
\newcommand{\sscap}[2]{#1^{}_{\scaleto{#2}{4pt}}}
\newcommand{\subfig}[1]{\textbf{#1}}

\title{\LARGE \bf
Optimal Path-Planning with 
Random Breakdowns
}

\author{Marissa Gee$^{1,2}$, Alexander Vladimirsky$^{2,3}$
\thanks{*Supported by the NSF DMS (awards 1645643, 1738010, and 2111522).}
\thanks{$^{1}$mag433@cornell.edu}
\thanks{$^{2}$Center for Applied Mathematics, Cornell University}%
\thanks{$^{3}$Department of Mathematics, Cornell University}%
}

\begin{document}
\maketitle
\thispagestyle{empty}
\pagestyle{empty}

\begin{abstract}
We propose a model for path-planning based on a single performance metric that accurately accounts for the the potential (spatially inhomogeneous) 
cost of breakdowns and repairs.
 These random breakdowns (or system faults) happen at a known, spatially inhomogeneous rate.
Our model includes breakdowns of two types: total, which halt all movement until an in-place repair is completed, and partial, after which movement continues in a damaged state
toward a repair depot.
We use the framework of piecewise-deterministic Markov processes to describe the optimal policy for all starting locations.
We also introduce an efficient numerical method that uses hybrid value-policy iterations to solve the resulting system of Hamilton-Jacobi-Bellman PDEs.
Our method is illustrated through a series of computational experiments that highlight the dependence of optimal policies on the rate and type of breakdowns, with one of them based on Martian terrain data near Jezero Crater.
\end{abstract}

\section{INTRODUCTION}
\label{sec:introduction}
Robustness is one of the central challenges in path-planning for autonomous vehicles.
Even when the dynamics and navigation environment are fully known, there is a question of how to account for the possibility of a partial or 
total 
breakdown. 
Much of the existing literature on planetary rovers models the undesirability of such events heuristically.  
The speed or cost functions might be modified ad hoc to reflect the risk~\cite{gennery1999traversability}, or terrain obstacles might be classified and essentially excluded from the environment before the actual trajectory planning~\cite{seraji1999traversability}.
A more rigorous approach could be based on multiobjective path-planning \cite{KumarVlad}, recovering the Pareto Frontier 
to capture all possible tradeoffs between the primary optimization criterion (e.g., the time to target when fully functional) and the risk of a breakdown along the way.
This, however, still ignores the possibility of repeated breakdowns.
More importantly, all of these approaches ignore that the actual 
consequences for the mission depend not just on the fact of a random breakdown but also on its {\em type} and {\em location}.
The primary contribution of this paper is a model that addresses these limitations rigorously and systematically in the framework of piecewise-deterministic Markov processes (PDMPs) \cite{davis1984piecewise}.

We consider an autonomous robot attempting to reach a target while subject to random breakdowns, whose (location-dependent) probability is known in advance.
After a {\em total breakdown}, the robot cannot move, and must pay for an in-place repair.
After a {\em partial breakdown}, it can continue traveling, but must pass through a repair depot on its way to the target.
Until a depot is reached, the damaged robot is essentially solving a different optimization problem since it might have a reduced speed or control authority, 
a restricted set of directions of motion, 
or a lower energy efficiency.
We model these different modes (fully functional and damaged) using PDMPs.
A PDMP is a stochastic model where at any point in time the system is in one of finitely many modes.
The system switches between these modes stochastically at known rates, while each mode specifies its own deterministic dynamics and running cost.
Recently, PDMPs have been applied to path planning problems with changing environments, where the modes represent environmental states such as changing wind direction 
\cite{ShenVlad, CarteeVlad_UQ}.
We take a similar approach here, but use the modes to model the status of the robot itself.

The range of path-planning methods used in the robotics and optimal control literature is truly broad (see \cite{Diks}\nocite{chazelle1985approximation}\nocite{kavraki1996probabilistic}\nocite{lavalle1998rapidly}-\cite{AltonMitchell3} for some examples).  Our approach is based on dynamic programming in continuous state and time: we obtain globally optimal trajectories by solving
Hamilton-Jacobi-Bellman PDEs.  This is a popular framework (e.g., \cite{garrido2016planning}\nocite{CarteeVlad_Poaching}-\cite{parkinson2020navigation}), particularly suitable when the state space is low-dimensional.  
But we note that our main ideas are also suitable for modifying  popular discrete state path-planning methods in higher dimensions~\cite{kavraki1996probabilistic, lavalle1998rapidly}.  


We limit our discussion to a simplified isotropic model of robot dynamics, primarily to streamline the exposition.
However, our framework is quite general and can be extended to the case of anisotropic dynamics or to more realistic, curvature-constrained models \cite{takei2013constrained}.
We start by describing the problem statement and the structure of the governing PDEs for three path-planning scenarios in Section \ref{sec:general_setting}.
We then present a novel iterative numerical method for solving these PDEs (Section \ref{sec:numerics}) and results of computational experiments (Section \ref{sec:numerical_results}).
Most of our test problems use synthetic data to illustrate the effects of system parameters on optimal trajectories, but the last example is more realistic and is based on terrain data 
for a region of Mars near Jezero crater~\cite{christensen2009jmars}.  We conclude by discussing future extensions in Section \ref{sec:conclusion}.

\section{GENERAL SETTING}
\label{sec:general_setting}

\subsection{Classical Path Planning}
We consider a robot that obeys the isotropic dynamics
\begin{equation}
\dot{\vy}(s) = \va(s) f(\vy(s)), \, \vy(0) = \vx,
\end{equation}
throughout some bounded domain $\Omega$, where $\va(s): \R \to S^1$ is the chosen direction of motion and $f$ is the speed of travel.
Henceforth, we use $\vx$ to refer to a generic point in $\Omega$ or a trajectory's initial condition and $\vy(s)$ to encode the dependence of the robot's position on time $s$.
We seek the policy $\va(\cdot)$ that minimizes 
the cumulative cost
\begin{equation}
J(\vx, \va(\cdot)) = \int_{0}^{T} K(\vy(s)) ds + q(\vy(T)).
\label{eq:classcost}
\end{equation}
Here $K$ is a running cost 
and $q$ is a terminal cost.
While $T$ could be set a priori, 
our focus is on
exit-time problems: let $G \subset \Omega$ be a finite set of points representing the target, then $T = \inf\{s|\vy(s) \in G\}$. 

The value function $u(\vx)$ is defined to encode the optimal cost-to-go from each point in $\Omega$:
\begin{equation}
u(\vx) = \inf_{\va(\cdot) \in \A} J(\vx, \va(\cdot)),
\end{equation}
where $\A$ is the set of measurable functions from $\R$ to $S^1$. 
Classical arguments from control theory \cite{CranLion} show that $u$ must be the viscosity solution of the Eikonal PDE (suppressing the dependence on $\vx$ for clarity):
\begin{equation}
|\nabla u| f = K, \,\, \vx \in \Omega \qquad u = q, \,\, \vx \in G,
\end{equation}
for which there are several well-known efficient numerical solvers; e.g., \cite{SethFastMarcLeveSet}\nocite{TsaiChengOsherZhao}\nocite{ChacVlad1}-\cite{ChacVlad2}.
Once $u$ is computed, the optimal policy can be recovered by setting $\va^*(s) = -\nabla u(\vy(s))/|\nabla u(\vy(s))|$.
In the special case that $K = 1$ and $q = 0$, the control $\va^*(s)$ is time-optimal.

\subsection{Simplified Model: Total Breakdowns Only}
Our first extension introduces the possibility of total breakdowns. 
We assume that breakdowns happen instantaneously and at random times. 
We view the time until the next breakdown as an exponential random variable with rate $\lambda >0$, except we allow $\lambda(\vx)$ to vary in space, making the risk of breakdown trajectory-dependent.
Let $R(\vx)$ be the in-place repair cost the robot must pay before continuing.
The \emph{expected cost} associated with a trajectory $\vy(s)$ is then given by
\begin{IEEEeqnarray}{rCl}
J(\vx, \va(\cdot)) &=& \int_{0}^{T} K(\vy(s)) + \lambda(\vy(s)) R(\vy(s)) ds
\label{eq:totalbdcost}
\end{IEEEeqnarray}
(where we now assume $q = 0$ at the target). 
While we model $R$ as being paid instantaneously, we will see later that we can choose $R$ to capture the time taken for a repair (assuming the robot does not accrue running cost while broken down).
In that case, $s$ becomes the time spent \emph{moving}, not the total time \emph{taken}.
The structure of \eqref{eq:totalbdcost} is identical to that of \eqref{eq:classcost}, and thus the same arguments show that the value function $u$ solves
\begin{equation}
|\nabla u| f = K + \lambda R, \,\, \vx \in \Omega \qquad u = 0, \, \vx \in G.
\end{equation}

We now describe a particular model for $R$.
Suppose the domain contains a finite number of repair depots at $D = \{\tilde{\vx}^{}_1, ..., \tilde{\vx}^{}_M \} \subset \Omega$, and that a repair vehicle must travel from a depot to the broken down robot to repair it. 
If the repair vehicle has its own speed $f_R(\vx)$ and running cost $K_R(\vx)$ and minimizes its own travel cost, we can write
\begin{gather}
R = \sscap{u}{R} + \sscap{R}{F}, \\
|\nabla \sscap{u}{R} | \sscap{f}{R} = \sscap{K}{R},\,\, \vx \in \Omega  \qquad \sscap{u}{R} = \sscap{R}{D}, \, \vx \in D \label{eq:uRpde}
\end{gather}
where $\sscap{R}{D}(\vx)$ and $\sscap{R}{F}(\vx)$ are depot and breakdown location-dependent repair costs, respectively.
This structure provides great flexibility in modeling the cost of a total breakdown, while still allowing $R$ to be precomputed throughout the domain.
For example, if we set $\sscap{K}{R} = 1$, $\sscap{R}{F}$ to be the time for an in-place repair, and $\sscap{R}{D}$ to be the pre-dispatch waiting time at each depot, then $R$ captures the total time spent broken down.
Alternatively, if the repair vehicle charges a rate $C$ for its time, we can scale $\sscap{R}{F}$ by $C$, $\sscap{K}{R}$ by $2C$ (to account for the return trip), and set $R_D = 0$ (assuming we are not charged for waiting time) to model the amount \emph{paid} for a repair.

\subsection{Full Model: Total and Partial Breakdowns}
Now, we allow for two types of breakdown: 
total, as defined above, and partial, where the robot can continue moving after the breakdown (possibly with 
a lower $f$ or higher $K$ or $\lambda$).
As before, a total breakdown immobilizes the robot until it is fixed in place.  Thus, we now have two distinct modes
in which the robot can travel and must decide on a trajectory:
mode 1, when it is fully functional, and mode 2, when it is damaged
and moving toward a repair depot.
We model this as a PDMP, where each mode has its own value function and they are coupled due to switching via breakdowns and repairs.
If we first assume $u_2$ (the value function in mode 2) is known, then the expected cost in mode 1 is
\begin{IEEEeqnarray}{rCl}
J_1(\vx, \va(\cdot)) &=& \E\left[\int_{0}^{T_1} K_1(\vy(s)) + \lambda_{1}(\vy(s)) R(\vy(s)) ds \nonumber \right.\\
&&+ \mathds{1}_{[T^{}_b < T^{}_{G}]} u_2(\vy(T_b))\Bigg]. 
\label{eq:cost1}
\end{IEEEeqnarray}
where $T_b$ is the time at which the first partial breakdown occurs, $\sscap{T}{G} = \inf\{s | \vy(s) \in G\}$ is the time the robot would reach the target if no breakdown occurred, and $T_1 = \min\{T_b, \sscap{T}{G}\}$.
The last term in \eqref{eq:cost1} is only nonzero when a partial breakdown occurs before the robot reaches the goal, 
in which case the robot switches to mode 2. 
If partial breakdowns occur at some known rate $\phi(\vx)$, then this is an example of a randomly-terminated finite-horizon control problem, as outlined in \cite{AndrewsVlad}.

We can similarly define $J_2$, assuming that $u_1$ is known. 
In mode 2, the robot returns to mode 1 by being repaired at a depot or having a total breakdown and paying $R$. 
Let $\sscap{T}{B}$ be the time of the first total breakdown, $\sscap{T}{D} = \inf\{s | \vy(s) \in D\}$ be the time the robot would reach the depot if no breakdowns occurred, and $T_2 = \min\{\sscap{T}{B}, \sscap{T}{D}\}$, then
\begin{IEEEeqnarray}{rCl}
J_2(\vx, \va(\cdot)) &=& \E\left[\int_{0}^{T_2} K_2(\vy(s)) ds \right.\nonumber\\
&& + \mathds{1}_{[T^{}_B < T^{}_D]} (R(\vy(\sscap{T}{B})) + u_1(\vy(\sscap{T}{B}))\\
&& + \mathds{1}_{[T^{}_B \geq T^{}_D]} (R_D(\vy(\sscap{T}{D})) + u_1(\vy(\sscap{T}{D}))\Bigg]\nonumber.
\end{IEEEeqnarray}
Here the last two terms again 
encode 
the optimal cost-to-go whenever a mode switch occurs. 
Following the derivation presented in \cite{CarteeVlad_UQ}, we arrive at the following system of coupled PDEs for the value functions:
\begin{IEEEeqnarray}{rClL}
|\nabla u_1| f_1 &=& K_1 + \lambda_{1} R + \phi (u_2 - u_1)\hspace{0.2cm}& \vx \in \Omega \label{eq:u1pde}\\
u_1 &=& 0 &\vx \in G \nonumber \\
|\nabla u_2| f_2 &=& K_2 + \lambda_{2} (R + u_1 - u_2)& \vx \in \Omega \label{eq:u2pde}\\
u_2 &=& \sscap{R}{D} + u_1 &\vx \in D. \nonumber 
\end{IEEEeqnarray}

The inclusion of partial breakdowns significantly complicates the model.
Previously it was possible to sequentially solve for $R$ and then $u$, but due to the coupling of $u_1$ and $u_2$ that strategy is no longer possible in general.
One exception is when the only depot is located at the target and there are no total breakdowns in mode 2 ($\lambda_2 = 0$).
In this case, the boundary condition for $u_2$ is known, since 
plugging in $u_1 = 0$ on $G$ gives $u_2 = \sscap{R}{D}$ on $D$, and thus the PDEs become decoupled.
We can also recover simpler models as special cases of 
\eqref{eq:u1pde} and \eqref{eq:u2pde}.
First, setting $\lambda_1 = \lambda_2 = \phi = 0$ recovers the classical path-planning problem without breakdowns.
If we instead set $\lambda_2 = \phi = 0$, we obtain the simplified model from the previous subsection.
Finally, setting $\lambda_1 = \lambda_2 = 0$ corresponds to 
an additional model with only partial breakdowns, which we examine in Examples 3 and 4.


\section{NUMERICS}
\label{sec:numerics}
\subsection{Discretized Equations}
The main computational challenge in solving \eqref{eq:u1pde} and \eqref{eq:u2pde} is the coupling between $u_1$ and $u_2$: without it, the PDEs could be solved entirely using existing numerical techniques. 
Still, we will make use of multiple existing methods, so we start by reviewing them and later address solving the coupled systems.
We will use the discretization $(x,y) \approx (x^{i}, y^{j})$, where we assume for convenience that the distances between gridpoints, $\Delta x$ and $\Delta y$, are constant.

When solving Eikonal PDEs, we use the Fast Marching Method (FMM) \cite{SethFastMarcLeveSet}, which relies on having a causal discretization of the PDE. 
We use one-sided finite difference approximations of the derivatives, given by
\begin{equation}
D^{ij}_{\pm x} u = \frac{\pm u^{i\pm1,j} \mp u^{ij}}{\Delta x}, \qquad D^{ij}_{\pm y} u = \frac{\pm u^{i,j\pm1} \mp u^{ij}}{\Delta y},
\end{equation}
to define the upwind difference operator in the $x$-direction
\begin{equation}
\D^{ij}_x u  = \min\{D^{ij}_{+x}u, -D^{ij}_{-x}u, 0 \}
\end{equation}
and its equivalent for the $y$-direction.
Thus, for the Eikonal equation, the discretized $u$ must satisfy
\begin{align}
 \sqrt{\left(\D^{ij}_{x}u\right)^2 + \left(\D^{ij}_{y} u\right)^2} &= \frac{K^{ij}}{f^{ij}}.
\label{eq:eikdis}
\end{align}
at each point in the discretized domain.
We can solve for $u$ by propagating the boundary values as outlined in \cite{SethFastMarcLeveSet}.

We also require solvers for a class of uncertain horizon problems described in \cite{AndrewsVlad}, which extends the FMM  to problems with random termination and \emph{known} terminal cost.
Equations \eqref{eq:u1pde} and \eqref{eq:u2pde} fit this framework if we assume that either $u_2$ or $u_1$, respectively, are known.
If that were the case, we would arrive at the following discretized equations: 
\begin{IEEEeqnarray}{l}
\sqrt{\left(\D^{ij}_{x}u_1\right)^2 + \left(\D^{ij}_{y} u_1\right)^2} \label{eq:u1dis}\\
\qquad \qquad = \frac{K_1^{ij} + \lambda_1^{ij}R^{ij}}{f_1^{ij}} - \frac{\phi^{ij}}{f_1^{ij}}(u_1^{ij} - u_2^{ij}) \nonumber \\
\sqrt{\left(\D^{ij}_{x}u_2\right)^2 + \left(\D^{ij}_{y} u_2\right)^2}\label{eq:u2dis}\\
\qquad \qquad = \frac{K_2^{ij}}{f_2^{ij}} - \frac{\lambda_2^{ij}}{f_2^{ij}}(u_2^{ij} - (u_1^{ij} + R^{ij})) \nonumber
\end{IEEEeqnarray}
and each value function could be found using an existing causal method, assuming the other were known.

\begin{algorithm}[t]
\SetKw{Compute}{Compute}
\SetKw{Solve}{Solve}
\SetKw{Initialize}{Initialize}
\SetKw{Set}{Set}

\KwIn{$tol$, termination threshold}
\KwIn{$\rho < 1$, policy evaluation threshold}

\Solve{Eq. \eqref{eq:uRpde} for $\sscap{u}{R}$ using FMM} \\
\Set{$R = \sscap{u}{R} + \sscap{R}{F}$}

\Set{$n=0$, $d^{(0)} = tol + 1$, $\delta = d^{(n)}$} \\
\Initialize{$u^{(0)}_1$ and $u^{(0)}_2$ with overestimates}\\

\While{$d^{(n)} > tol$}{
\While{$d^{(n)} > \rho \cdot \delta$} {
\Initialize{$u^{(n+1)}_2 = u^{(n)}_1$ on $D$} \\
\Solve{Eq. \eqref{eq:u2pde} for $u^{(n+1)}_2$ using modified FMM with $u_1 = u_1^{(n)}$} \\
\Solve{Eq. \eqref{eq:u1pde} for $u^{(n+1)}_1$ using modified FMM with $u_2 = u_2^{(n+1)}$} \\
\Set{$n = n+1$, $d_n = \|u^{(n)}_1 - u^{(n-1)}_1 \|_\infty$} 
}
\Solve{Eqs. \eqref{eq:r1pde} and \eqref{eq:r2pde} for $r_1$ and $r_2$}\\
\Set{$\delta = \|u^{(n)}_1 - r_1 \|_\infty$, $u^{(n)}_1 = r_1$, $u^{(n)}_2 = r_2$} \\

}
\caption{Solve equations \eqref{eq:u1pde} and \eqref{eq:u2pde} via value-policy iteration.}
\label{alg:polvalit}
\end{algorithm}

\subsection{Value Iterations}
To take advantage of the methods presented above, we take an iterative approach similar to standard value iteration.
At each iteration, we first freeze the value of $u_2$ and use an existing method to solve \eqref{eq:u1pde} for $u_1$.
We then use the new version of $u_1$ to solve \eqref{eq:u2pde} for $u_2$ in the same way.
We alternate updates in this manner until the change between iterations, $\delta$, falls below a specified tolerance $tol$. 
The convergence can be proved by interpreting this algorithm as a (mode-by-mode) Gauss-Seidel relaxation of standard value iterations \cite{Bertsekas_DPbook}.

To start the iterative process, we need to initialize $u_2$ using some overestimates,
which we obtain by posing simpler ``pessimistic''  
problems. 
For all pessimistic problems, we assume that the robot does not leave mode 2 after a total breakdown, even after paying $R$ --- thus, the only coupling between $u_1$ and $u_2$ is on $D$.
This provides an overestimate if $f_2 \leq f_1$, $K_2 \geq K_1$, and $\lambda_2 \geq \lambda_1$.
To obtain a pessimistic version of $u_1$ at every $\tilde{\vx}^{}_m \in D$ (and thus an overestimate of the boundary conditions for \eqref{eq:u2pde}), we replace $f_k$, $\phi$, and $K_k + \lambda_k R$ by their worst values and assume that the vehicle moves along a straight line from $\tilde{\vx}^{}_m$ to the closest point in $G$, returning to $\tilde{\vx}^{}_m$ along the same line in case of a partial breakdown.  
This makes the problem essentially 1D, making it easy to solve analytically using a system of two coupled linear ODEs.
We then initialize $u_2$ throughout $\Omega$ by solving
\begin{equation}
|\nabla u_2| f_2 = K_2 + \lambda_2 R, \, \vx \in \Omega \qquad u_2 = R_D + \hat u_1, \, \vx \in D
\end{equation}
where $\hat u_1$ is the overestimate produced above.

\subsection{Acceleration via Value-Policy Iteration}
One shortcoming of value iterations is that $\delta$ quickly becomes small, even when the current value function is far from the solution.
Thus, pure value iterations are slow to converge, especially for poor initializations of the system.
To remedy this, we extend to PDMPs the method of combined value iterations and policy evaluations, outlined in \cite{grune2003programming} for single mode problems.
Policy evaluation is commonly found as a step in policy iteration, a popular dynamic programming technique, and involves fixing the control $\va(s)$ and computing the value function exactly for that fixed, suboptimal policy.

Fixing $\va_1(s)$ and $\va_2(s)$ in equations \eqref{eq:u1pde} and \eqref{eq:u2pde} we get
\begin{IEEEeqnarray}{rClL}
-(\va_1 \cdot \nabla r_1) f_1 &=& K_1 + \lambda_{1} R + \phi (r_2 - r_1)\hspace{0.2cm} &  \vx \in \Omega \label{eq:r1pde}\\
r_1 &=& 0 &\vx \in G \nonumber \\
-(\va_2 \cdot \nabla r_2) f_2 &=& K_2 + \lambda_{2} (R + r_1 - r_2)&\vx \in \Omega \label{eq:r2pde}\\
r_2 &=& R_D + r_1 &\vx \in D \nonumber, 
\end{IEEEeqnarray}
a system of coupled linear PDEs that can be discretized using finite differences and solved efficiently using any large-scale linear solver.
The combined value-policy iteration algorithm for two coupled value functions is outlined in Algorithm \ref{alg:polvalit}.

\section{NUMERICAL RESULTS}
\label{sec:numerical_results}

\begin{figure}
\begin{center}
  \begin{subfigure}{0.475\linewidth}
  \caption{}
  \label{fig:sub1}
  \includegraphics[width=\linewidth]{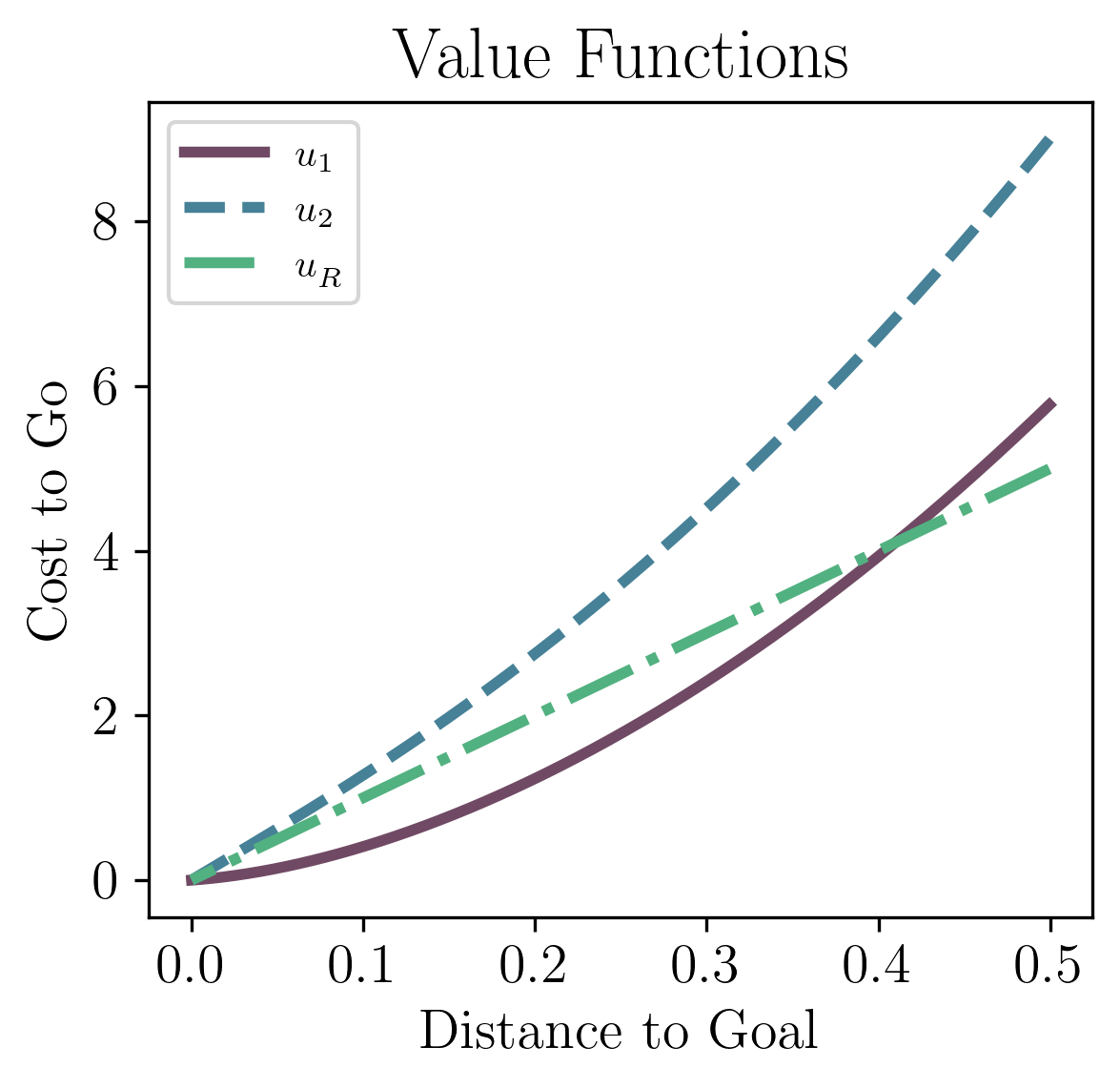}
  \end{subfigure}
  \begin{subfigure}{0.505\linewidth}
  \caption{}
  \label{fig:convb}
  \includegraphics[width=\linewidth]{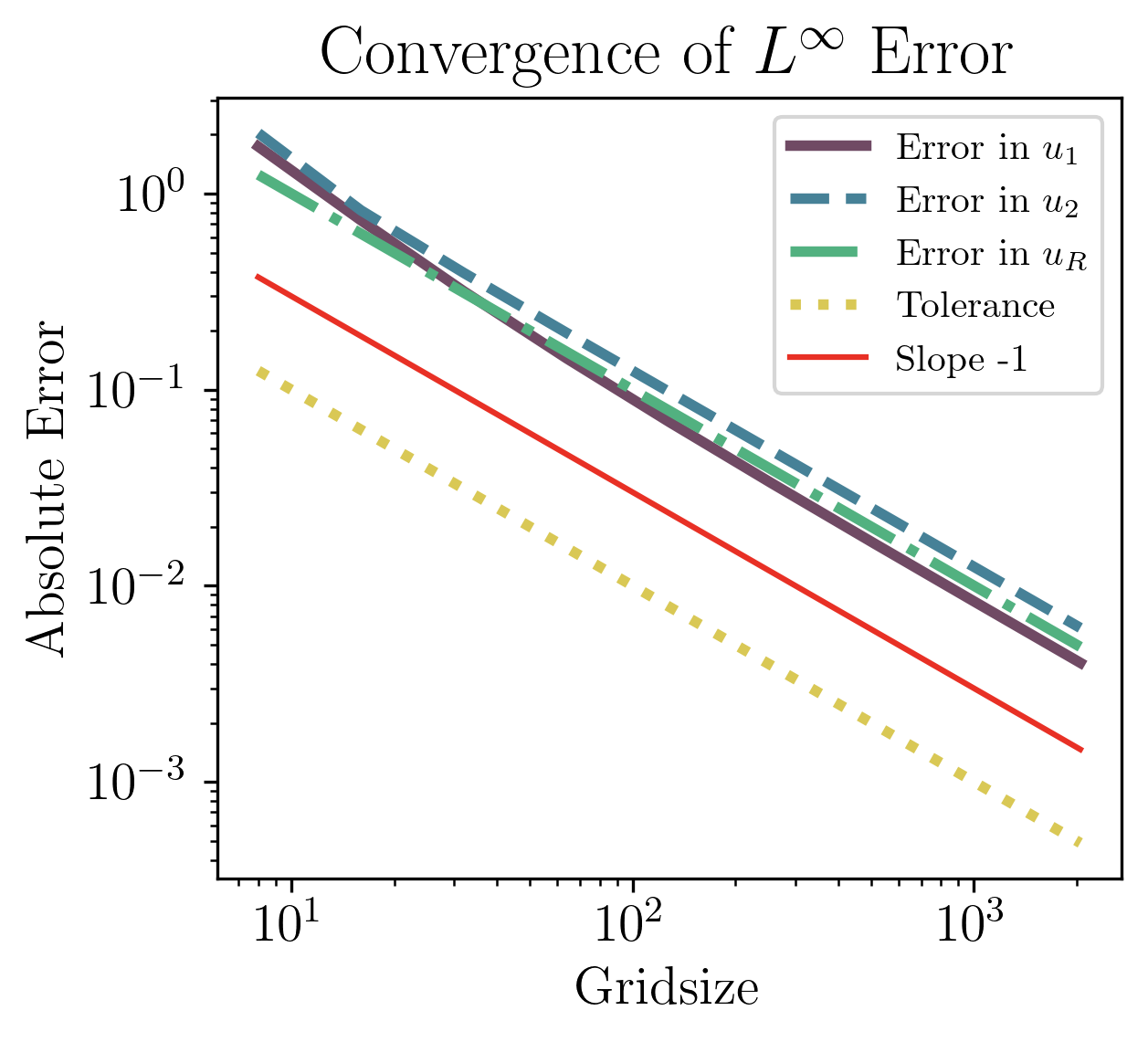}
  \end{subfigure}
\end{center}
\caption{Example 1. Radially symmetric value functions with $G = D = \{(0.5, 0.5)\}$. \textbf{(a)} Value functions for $\phi = 5$, $\lambda_1 = 0.5$, $\lambda_2 = 1.5$, $f_1 = 1$, $f_2 = 0.2$, and $f_R = 0.1$. \textbf{(b)} $L^\infty$ error between analytic and numeric solutions, with value iteration tolerance and slope $-1$ line for reference. Error exhibits first-order convergence.}
\label{fig:conv}
\end{figure}
 
For all examples except the first, the value functions are represented on a 501$\times$501 grid on $[0,1] \times [0,1]$.
For all examples, we set $tol = 2/(\Delta x + \Delta y)$, $K = 1$ in both modes, $R_F = 1$, $R_D = 0$, and use spatially homogenous environmental parameters $f_1$, $f_2$, $f_R$, $\lambda_1$, $\lambda_2$, and $\phi$ unless otherwise specified. \footnote{Code for all examples is available at \url{https://github.com/eikonal-equation/Random_Breakdowns}}

\subsection{Example 1: Convergence of Iterative Scheme}
We start with
a model for which the solution can be computed analytically: one target and one depot, both at the same location.
The 
2D solution is radially symmetric about $G$, and can be obtained by solving a system of ODEs analytically.
Figure \ref{fig:conv}\subfig{a} shows the value functions 
as functions of distance to the target.
Figure \ref{fig:conv}\subfig{b} shows the first-order convergence of the iterative scheme under grid refinement.

\begin{figure}
\begin{center}
  \begin{subfigure}{0.49\linewidth}
  \caption{}
  \label{fig:obstacleMBDphi}
  \includegraphics[width=\linewidth]{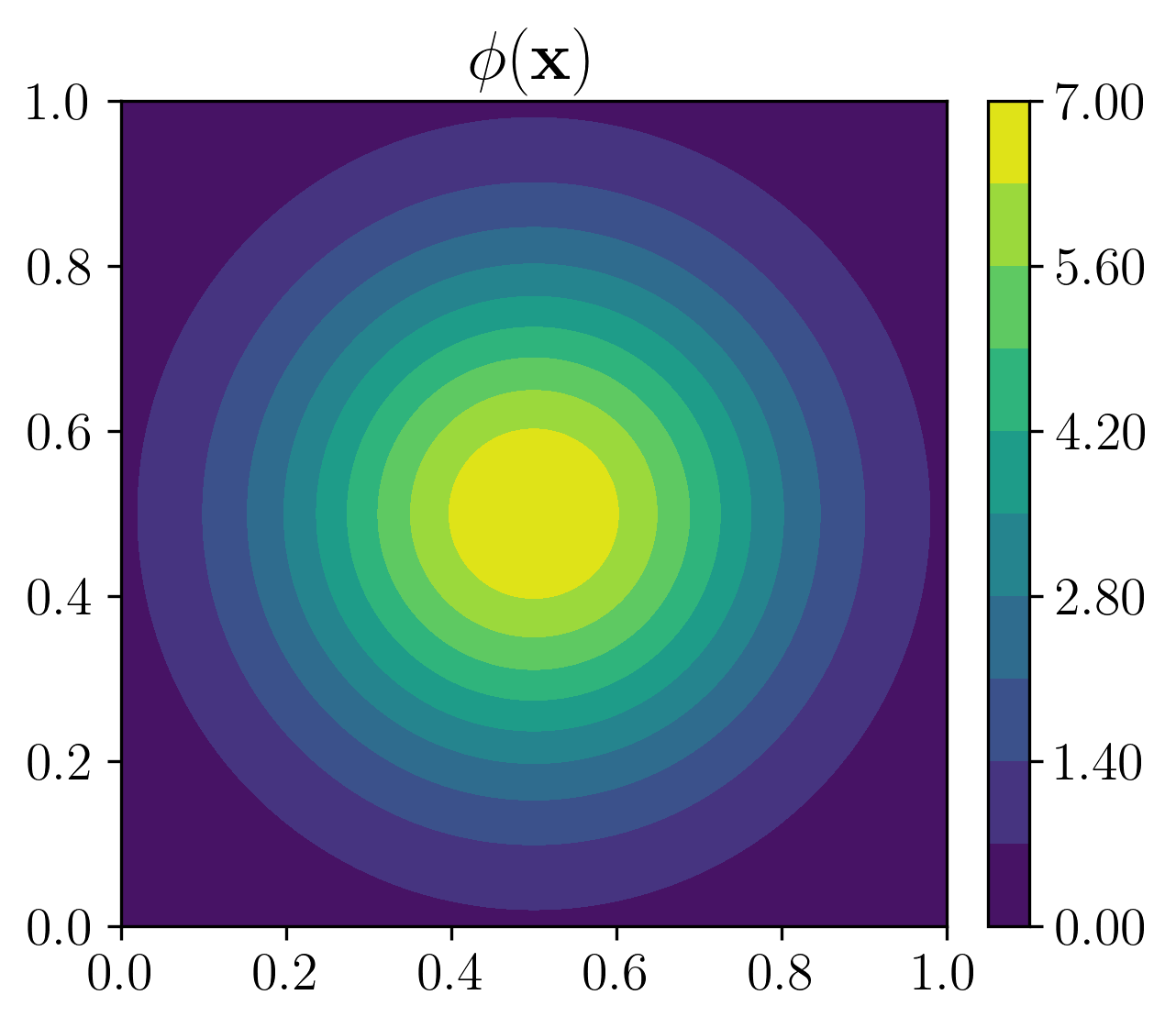}
  \end{subfigure}
  \begin{subfigure}{0.49\linewidth}
  \caption{}
  \label{fig:obstacleMBDu1NMB}
  \includegraphics[width=\linewidth]{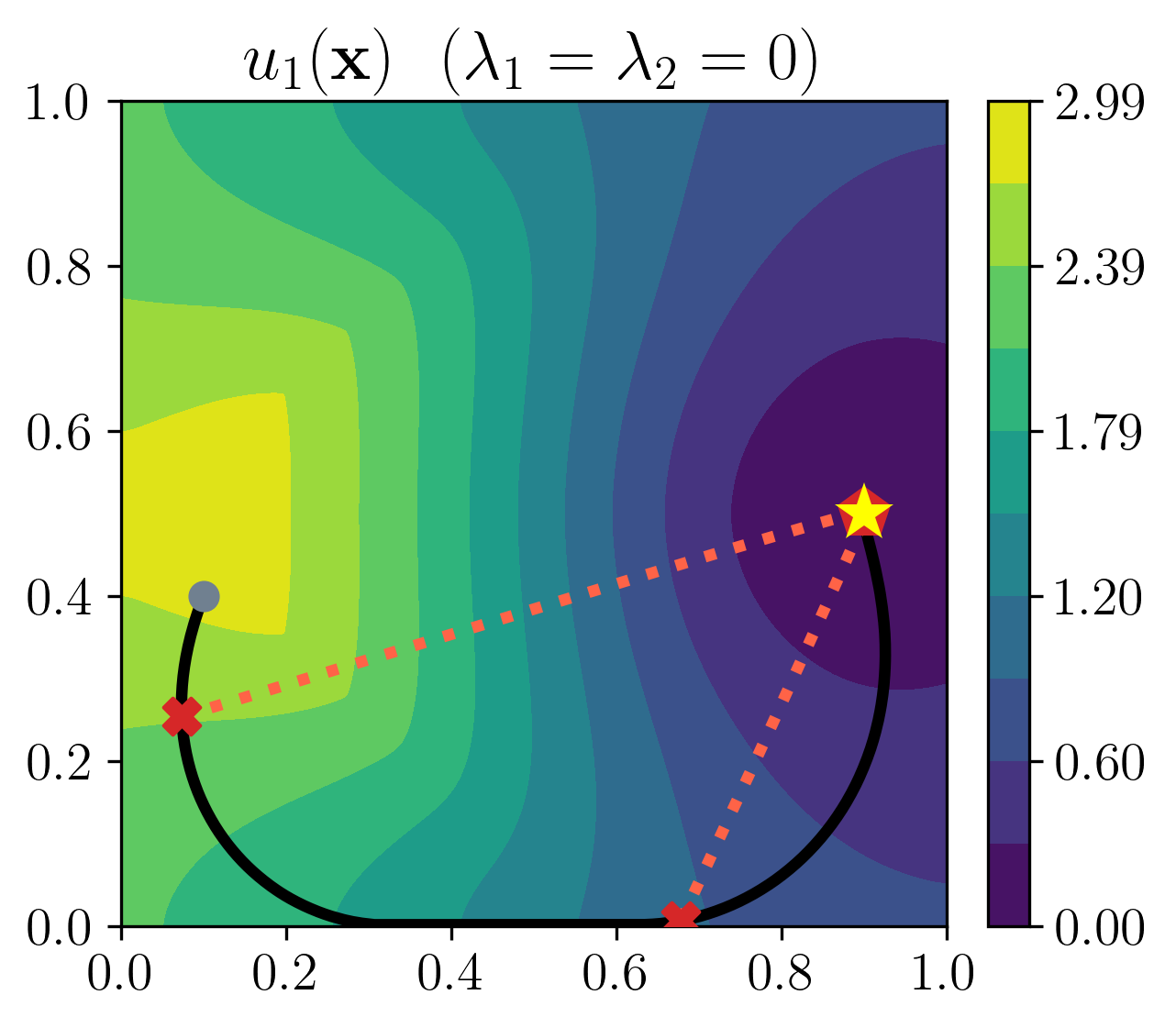}
  \end{subfigure}
  \begin{subfigure}{0.49\linewidth}
  \caption{}
  \label{fig:obstacleMBDu2}
  \includegraphics[width=\linewidth]{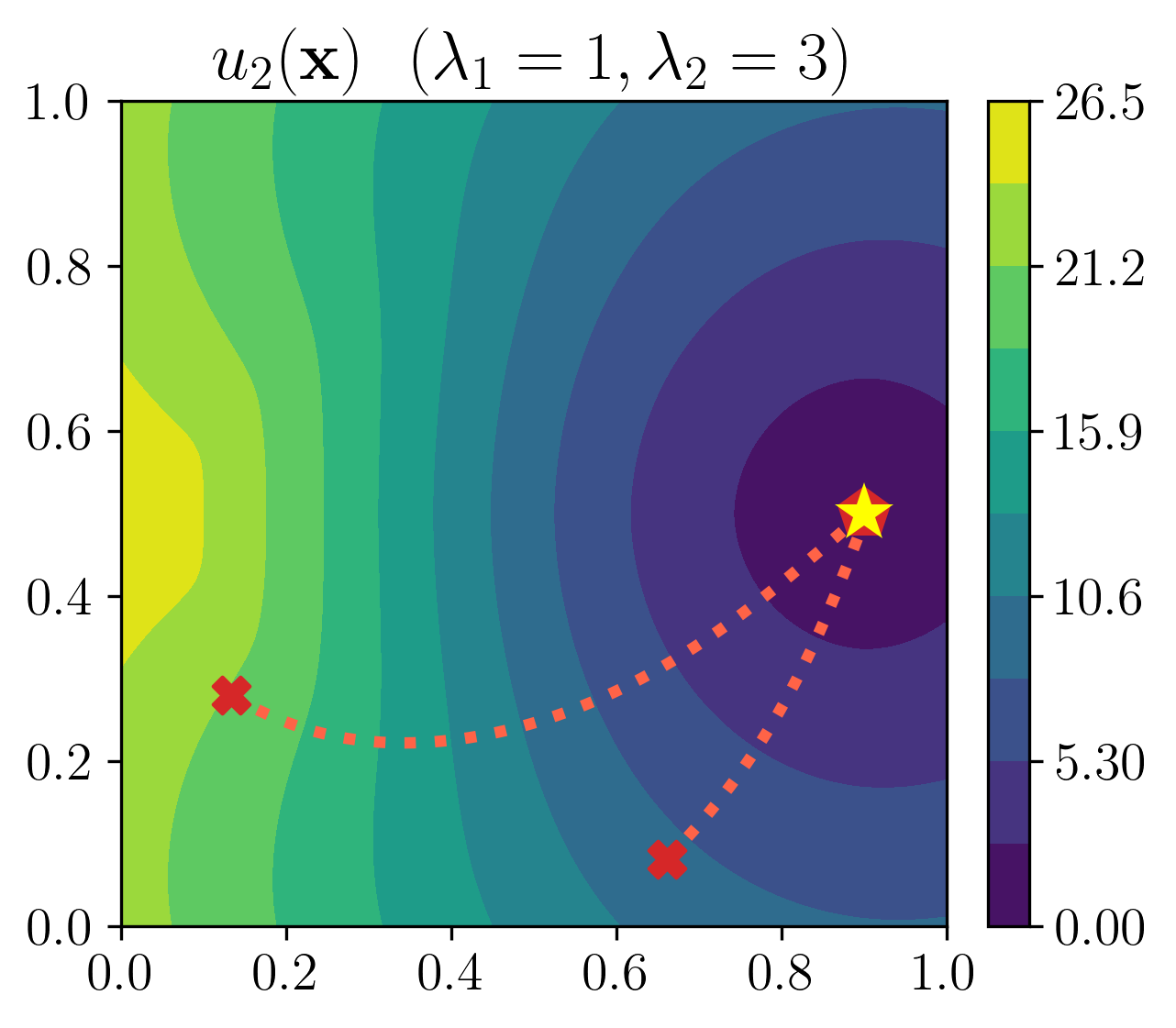}
  \end{subfigure}
  \begin{subfigure}{0.49\linewidth}
  \caption{}
  \label{fig:obstacleMBDu1}
  \includegraphics[width=\linewidth]{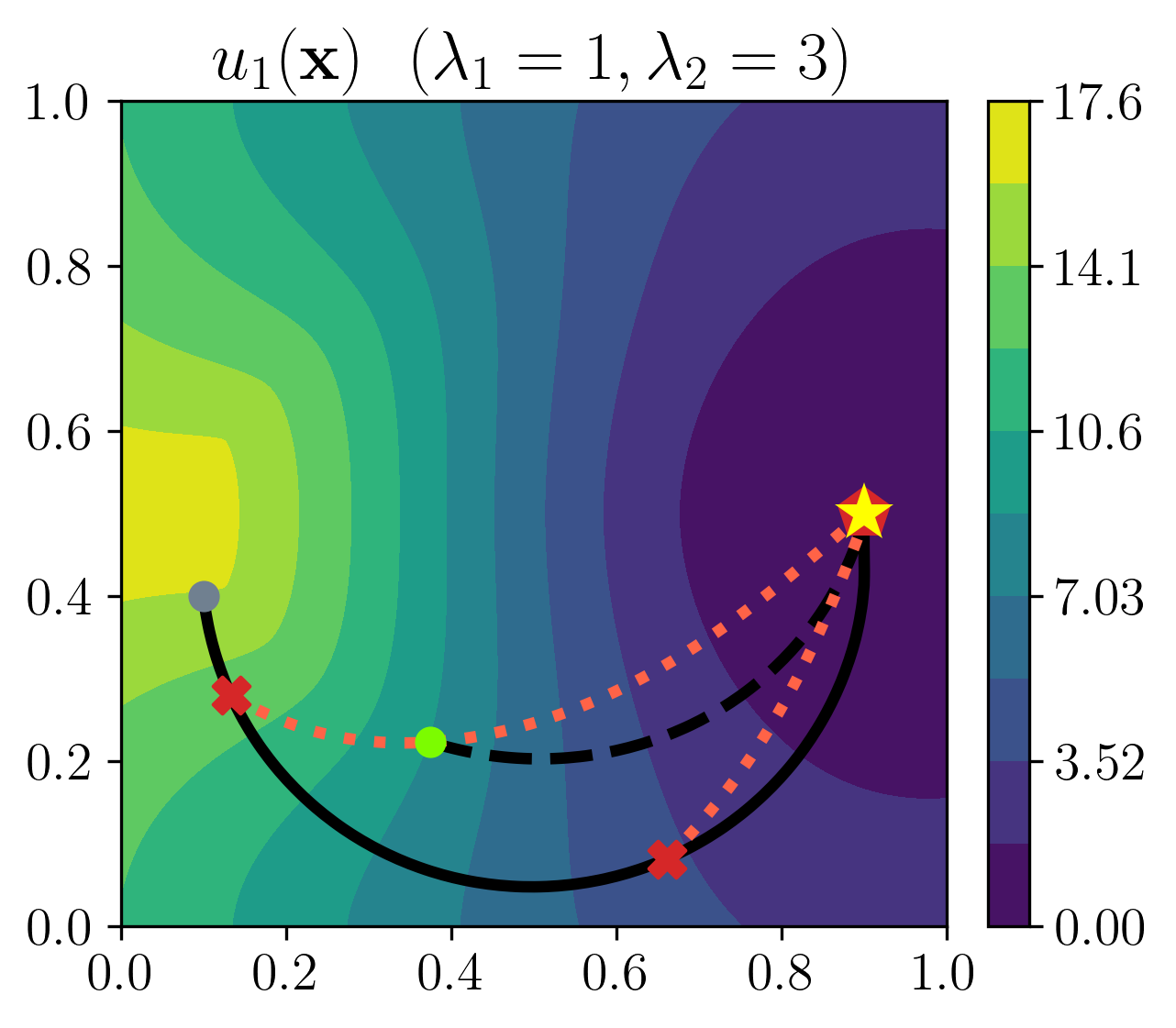}
  \end{subfigure}
\end{center}
\caption{For all figures: gold stars are targets, red pentagons are depots, grey dots are sample starting locations, red Xs are potential breakdown locations, and green dots are potential in-place repairs. Black trajectories are optimal in mode 1 and red dotted trajectories are optimal in mode 2. Example 2: $G = D = \{(0.9, 0.5)\}$. \textbf{(a)} The partial breakdown rate. \textbf{(b)} Value function for mode 1 and sample trajectories for modes 1 and 2 when $\lambda_1 = \lambda_2 = 0$. \textbf{(c)} Value function and sample trajectories for mode 2, when $\lambda_1 = 1$ and $\lambda_2 = 3$. The optimal trajectories are orthogonal to the level sets of $u_2$. \textbf{(d)} Value function and sample trajectories for mode 1 when $\lambda_1 = 1$ and $\lambda_2 = 3$.}
\label{fig:obstacleMBD}
\end{figure}

\subsection{Example 2: Inhomogenous Partial Breakdown Rate}
We now consider an inhomogenous partial breakdown rate, shown in Figure \ref{fig:obstacleMBD}\subfig{a}.
We set $f_1 = 1$, and use $f_2 = 0.2 f_1$ and $\sscap{f}{R} = 0.1 f_1$ to enforce slower speeds for the damaged robot and the repair vehicle.
Figure \ref{fig:obstacleMBD}\subfig{b} shows the mode 1 value function in the absence of total breakdowns. 
The optimal trajectories clearly avoid areas with high values of $\phi$.
Without total breakdowns ($\lambda_1 = \lambda_2 = 0$), there is no need to compute $u_R$ and $u_2$ is radially symmetric about $G$, so it is optimal to travel in a straight line to the target in mode 2.

This changes once we introduce a chance of total breakdowns.
Figures \ref{fig:obstacleMBD}\subfig{c} and \ref{fig:obstacleMBD}\subfig{d} show the value functions 
for a \emph{constant} chance of total breakdown ($\lambda_1 = 1$, $\lambda_2 = 3$). 
Now, the optimal trajectories in mode 2 are no longer straight lines.
Instead, the robot avoids areas where $\phi$ is high so as not to be caught there if it is repaired in-place after a total breakdown.
Thus, the introduction of even a constant value of $\lambda_2$ can disrupt the radial symmetry of $u_2$ about the target when $\phi$ is inhomogenous.

\subsection{Example 3: Changing Partial Breakdown Rate}
\begin{figure}
\begin{center}
\begin{subfigure}{0.49\linewidth}
	\caption{}
	\label{fig:changeBDRa}
	\includegraphics[width=\linewidth]{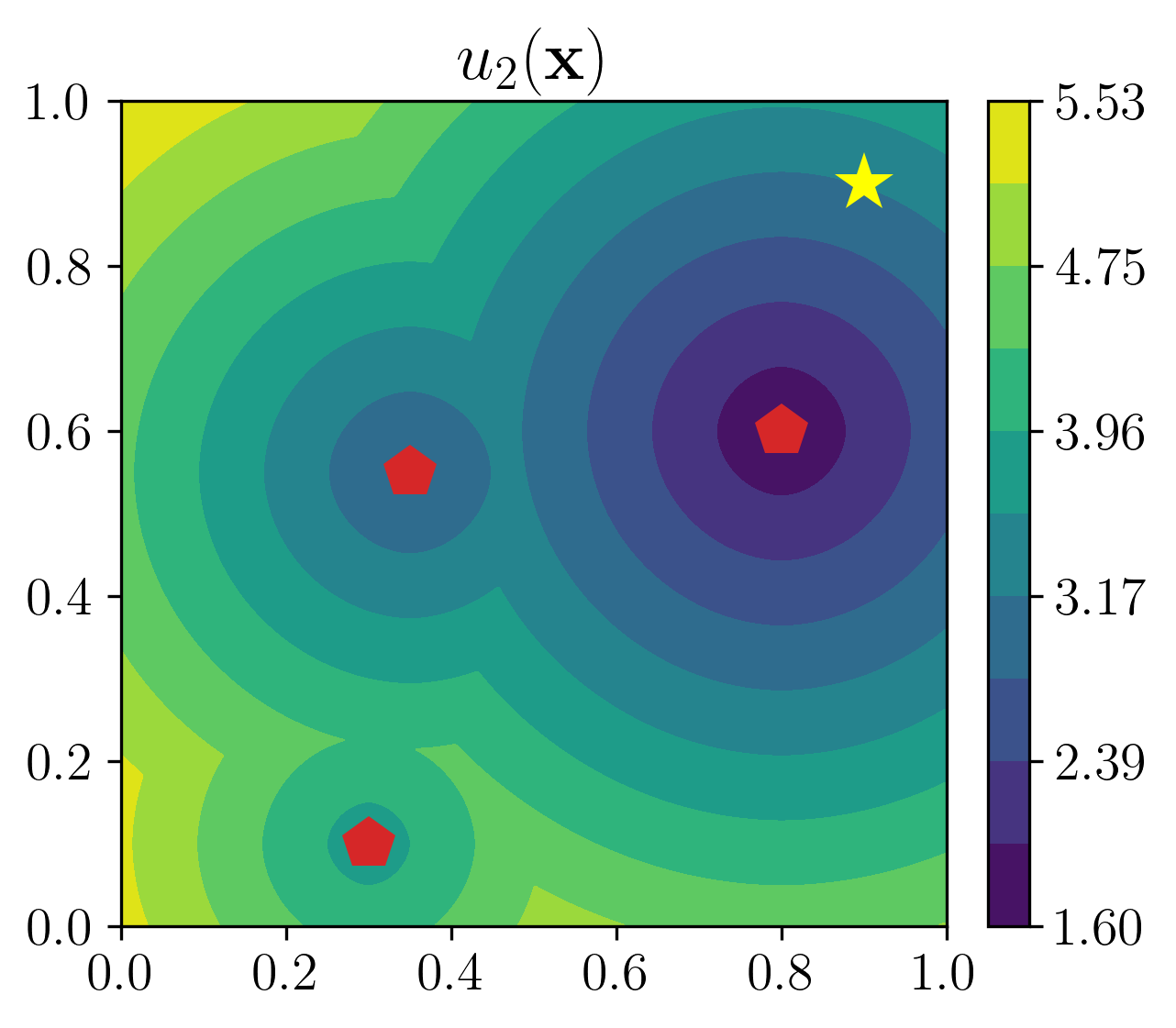}
\end{subfigure}
\begin{subfigure}{0.49\linewidth}
	\caption{}
	\label{fig:changeBDRb}
	\includegraphics[width=\linewidth]{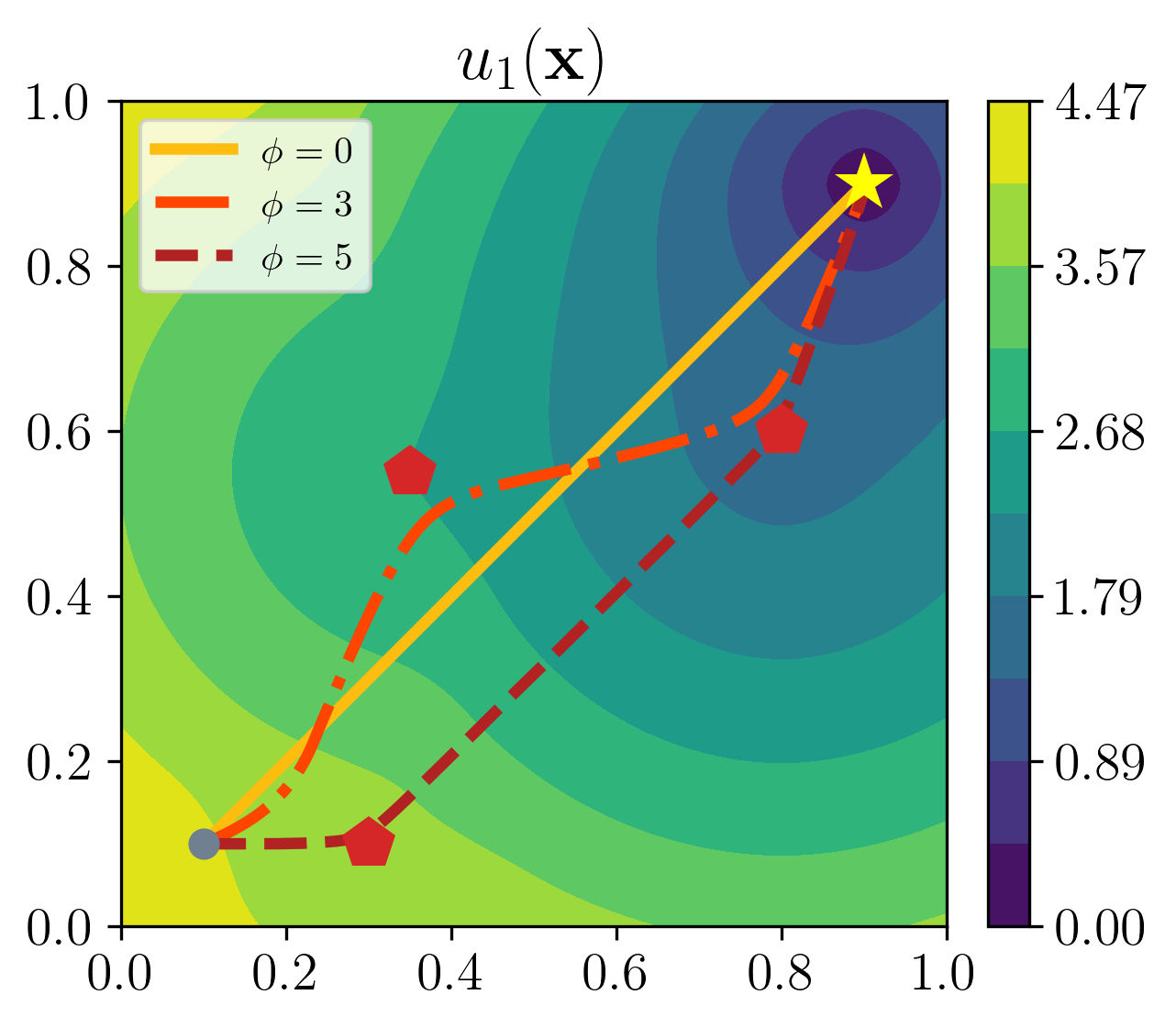}
\end{subfigure}
\end{center}
\caption{Example 3, $G = \{(0.9, 0.9\})$ and $D = \{(0.3,0.1), (0.35, 0.55), (0.8, 0.6)\}$. \textbf{(a)} Mode 2 value function with $\phi = 3$ and $\lambda_1 = \lambda_2 = 0$. \textbf{(b)} Mode 1 value function for the same parameters. Sample optimal trajectories  in mode 1 for $\phi = 0$ (yellow), $\phi = 3$ (orange \& orthogonal to the shown level sets of $u_1$), and $\phi=5$ (red) starting at $\vx=(0.1, 0.1)$.}
\label{fig:changeBDR}
\end{figure}

We now examine
how changes in $\phi$ 
affect the optimal paths.
We set $\lambda_1 = \lambda_2 = 0$, and use the same speed values as Example 2 but with three different depots (Figure \ref{fig:changeBDR}\subfig{a}).
Since $G \bigcap D = \emptyset,$ a damaged robot \emph{must} 
travel to a depot rather than directly to the target. 
This example illustrates the tradeoffs associated with extending the planned mode 1 trajectories.
Longer trajectories can remain closer to $D$, which 
reduces the expected time spent traveling in mode 2, but also increases the chance that a breakdown occurs at all. 

Figure \ref{fig:changeBDR}\subfig{b} shows the impact of changing $\phi$ on the optimal trajectories in mode 1.
With no chance of breakdowns ($\phi=0$), the robot would follow the yellow straight line to $G.$
For relatively low values of $\phi>0,$ it is optimal to extend the path (orange) in order to stay close to $D$ and decrease the time spent traveling in mode 2.
However, as $\phi$ increases, 
breakdowns are more frequent and the robot may have to travel in mode 2, even after a recent repair. 
Eventually, it is no longer worth extending the path (red) to be near every repair depot, 
and it becomes optimal to prioritize reaching the best depot (the one with the lowest expected cost-to-go) as quickly as possible.

\subsection{Example 4: Real-world Terrain}
\begin{figure*}
\begin{subfigure}{0.305\linewidth}
	\caption{}
	\label{fig:marsExamplef2}
	\includegraphics[width=\linewidth]{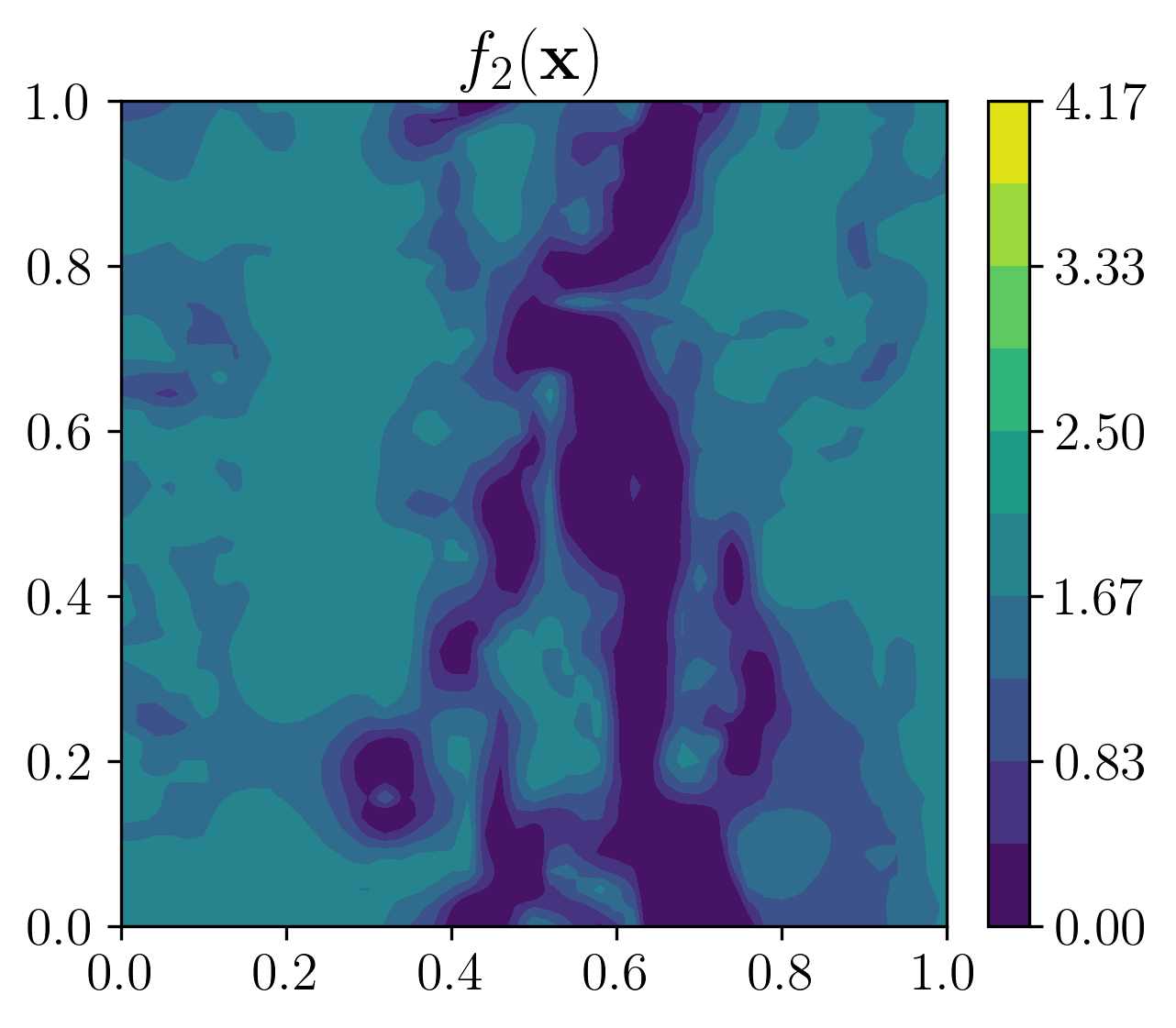}
\end{subfigure}
\begin{subfigure}{0.305\linewidth}
	\caption{}
	\label{fig:marsExamplef1}
	\includegraphics[width=\linewidth]{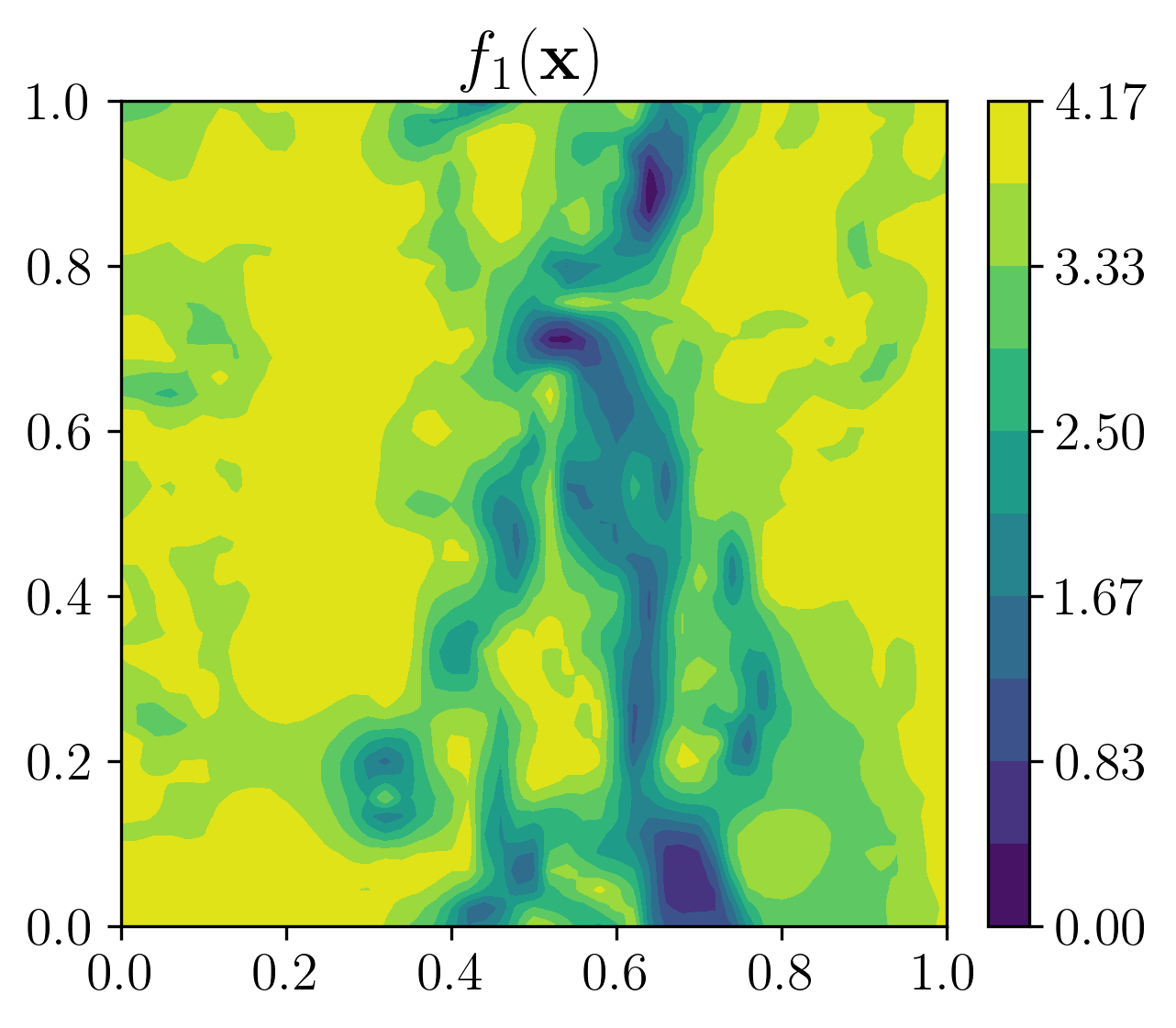}
\end{subfigure}
\begin{subfigure}{0.305\linewidth}
	\caption{}
	\label{fig:marsExamplephi}
	\includegraphics[width=\linewidth]{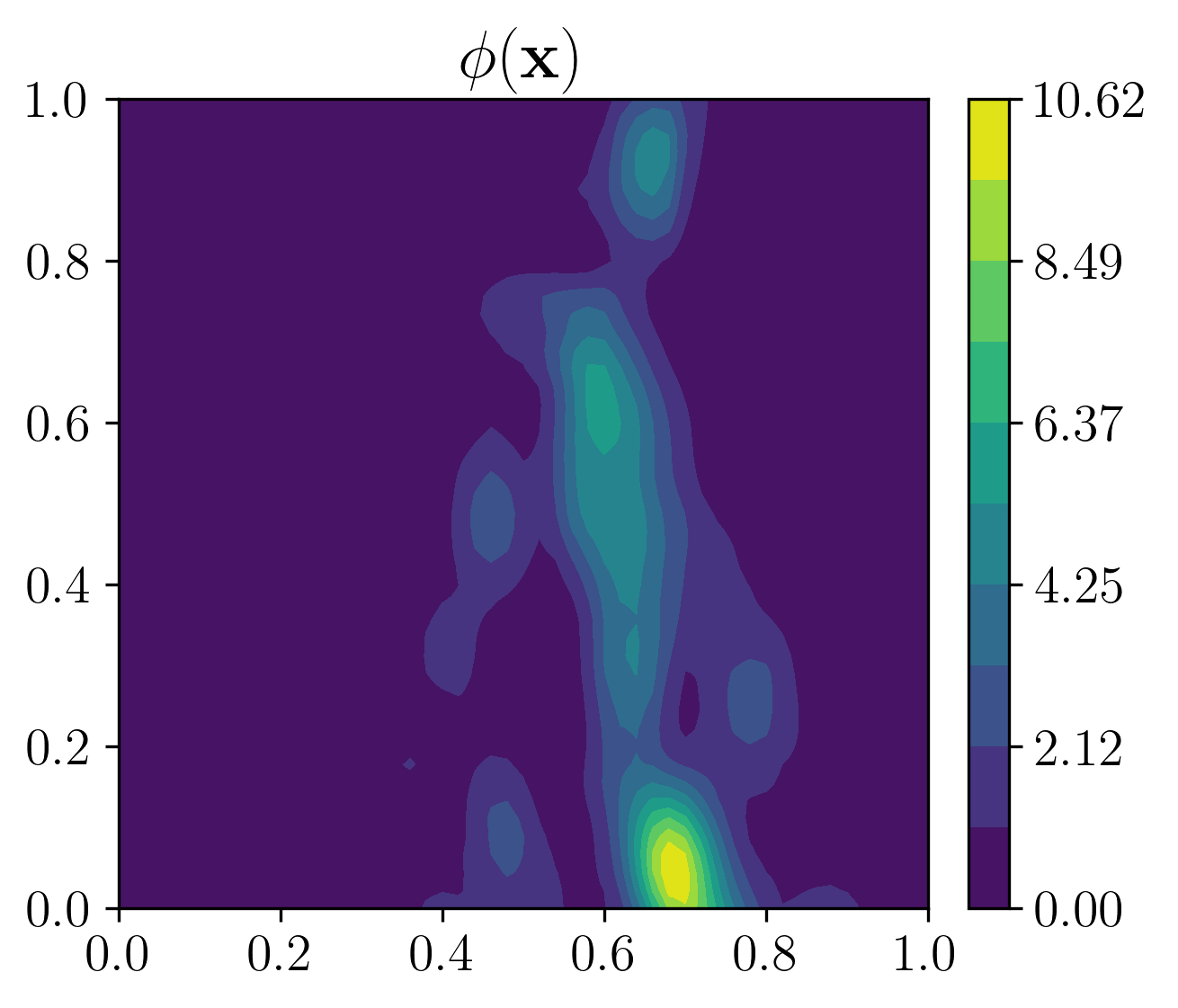}
\end{subfigure}\\
\begin{subfigure}{0.305\linewidth}
	\caption{}
	\label{fig:marsExampleu2}
	\includegraphics[width=\linewidth]{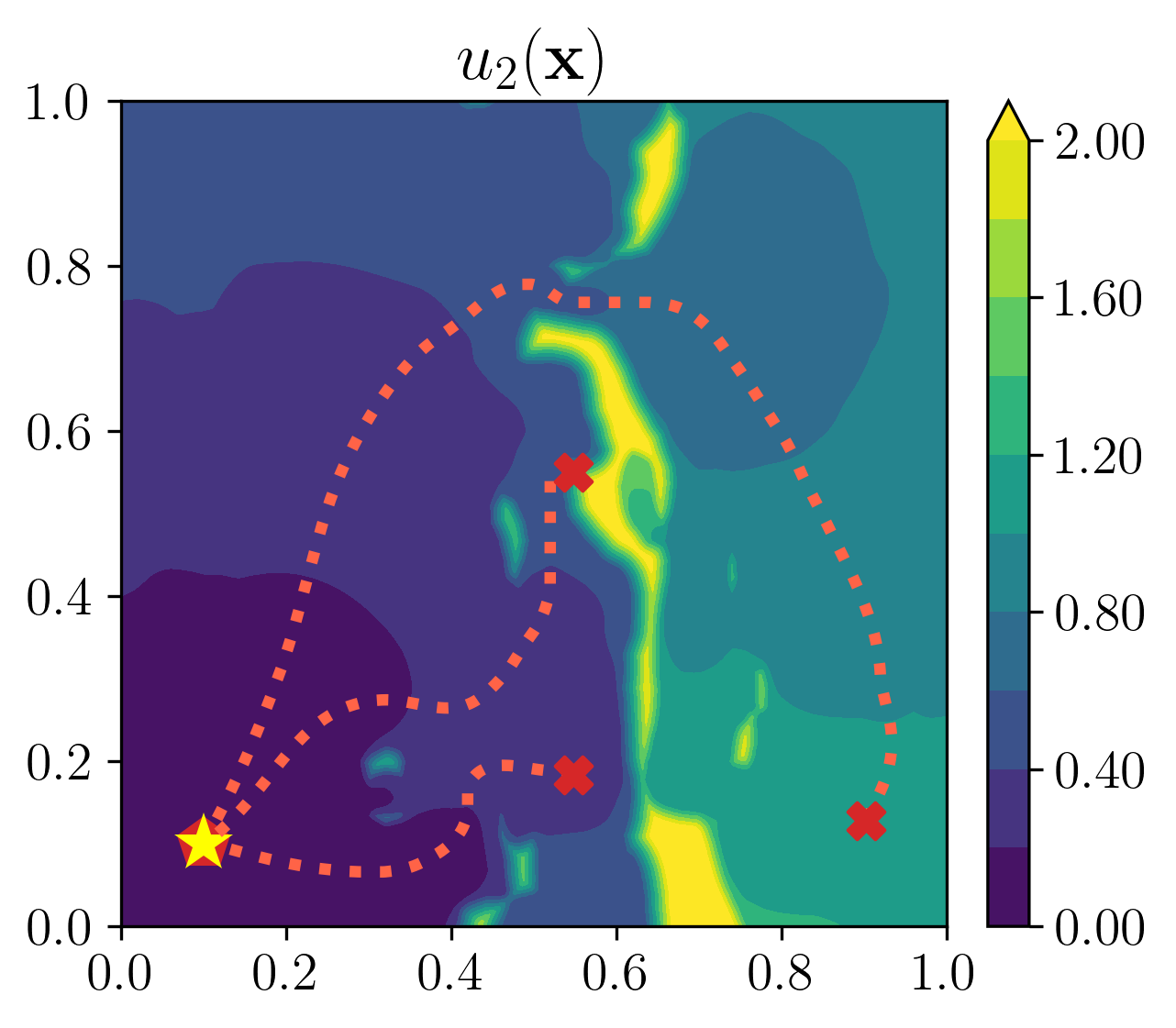}
\end{subfigure}
\begin{subfigure}{0.305\linewidth}
	\caption{}
	\label{fig:marsExampleu1}
	\includegraphics[width=\linewidth]{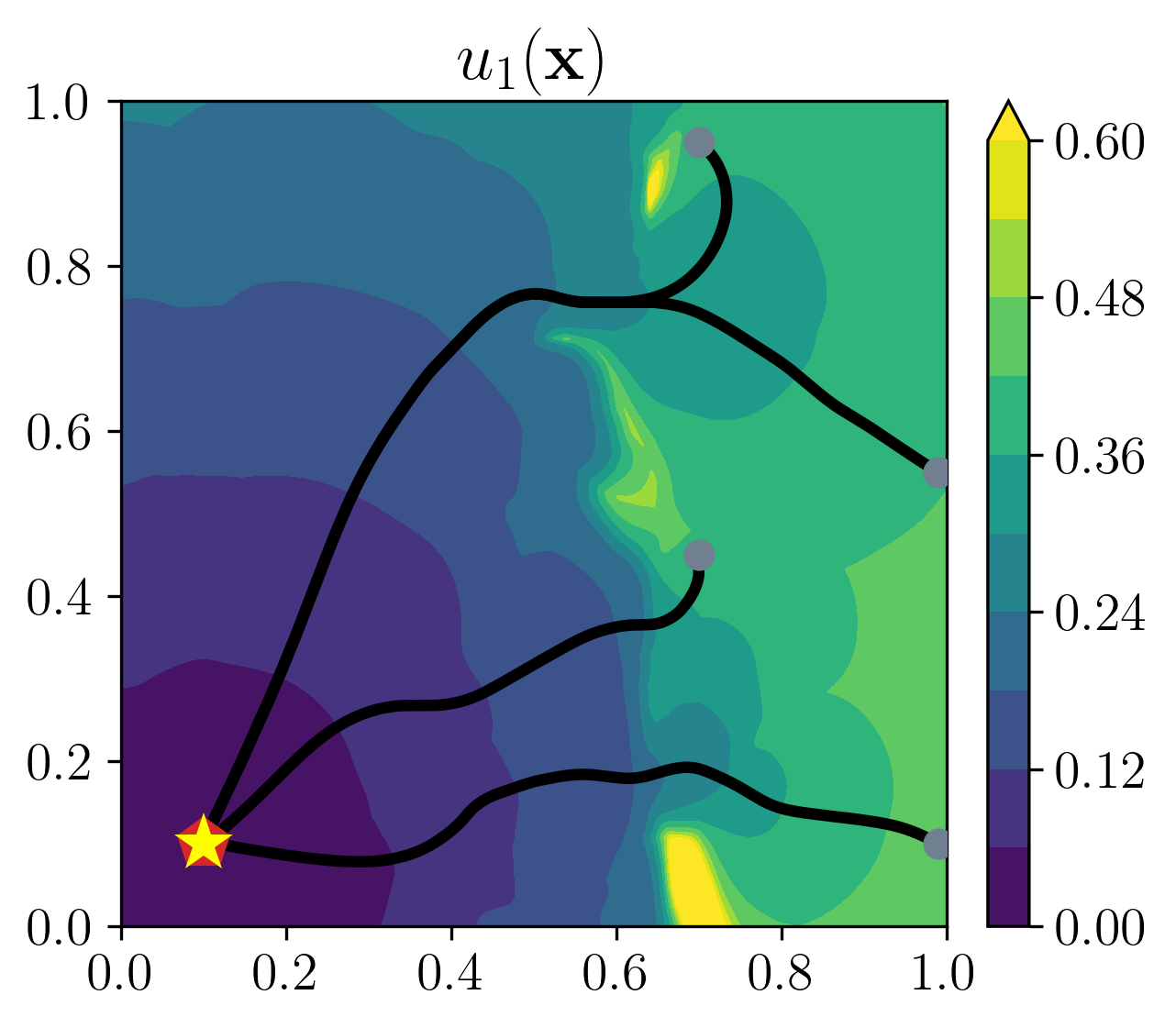}
\end{subfigure}
\begin{subfigure}{0.35\linewidth}
	\caption{}
	\label{fig:marsExampleEle}
	\includegraphics[width=\linewidth]{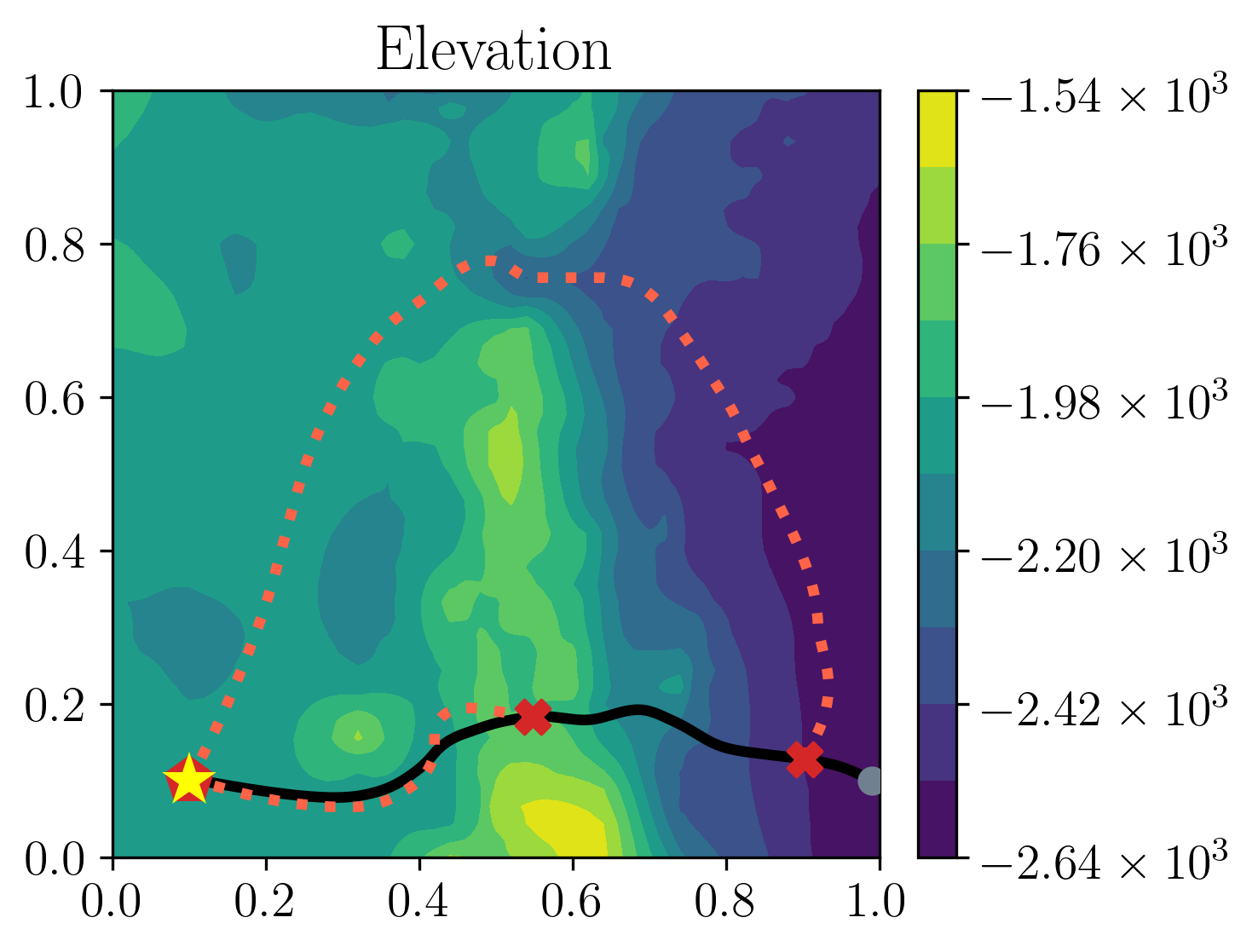}
\end{subfigure}
\caption{Example 4, $G = D = \{(0.1, 0.1)\}$. \textbf{(a)} Speed of rover in mode 2. \textbf{(b)} Speed of rover in mode 1. \textbf{(c)} Partial breakdown rate. \textbf{(d)} Value function and sample trajectories for mode 2. \textbf{(e)} Value function and sample trajectories for mode 1. \textbf{(f)} Elevation and sample optimal trajectories in modes 1 and 2.}
\label{fig:marsExample}
\end{figure*}

As our last example, we consider
path-planning for a Mars rover. 
We examine a region of Mars near Jezero crater, which was the landing site of NASA's Perseverance rover, and thus has extensively characterized terrain.
All data was accessed via JMARS, a Mars GIS \cite{christensen2009jmars}. 
We 
base the rover's speed on the slope of the terrain and measure it in units of meters/sol (a \emph{sol} is a Martian day).
Current Mars rovers can navigate slopes of up to approximately $\overline{\sigma} = 20^\circ$, and travel more slowly on steeper slopes due to increased wheel slippage 
\cite{heverly2013traverse}.
For simplicity we assume that the rover travels at a small but nonzero speed $f_{min} = 1\, m/sol$ for any slope above $\overline{\sigma}$, and following \cite{ono2016data} we set $f_{max} = 200\, m/sol$, which yields
\begin{align}
f_1(\vx) = \begin{cases}
f_{max}-\frac{f_{max} - f_{min}}{\overline{\sigma}}\sigma(\vx) & \sigma \leq \overline{\sigma} \\
f_{min} & \sigma > \overline{\sigma}
\end{cases}
\end{align}
where $\sigma(\vx)$ is the slope. 
We assume that $\overline{\sigma}$ and $f_{max}$ are decreased by a factor of $2$ in mode 2.
Figures \ref{fig:marsExample}\subfig{a} and \ref{fig:marsExample}\subfig{b} show $f_1$ and $f_2$ for the region we consider.

We do not allow total breakdowns ($\lambda_{1} = \lambda_2 = 0$), since repairs are not feasible for a Mars rover.
However, we do allow for partial breakdowns, and set $D = G$, ensuring that the rover will continue to travel to the target in mode 2.
The breakdown rate is modeled as proportional to the \emph{terrain roughness}, $\rho(\vx)$.
We compute $\rho$ as the root-mean-square height, a standard measure in geology \cite{shepard2001roughness} that has also been used to characterize traversability for Martian rovers in the past \cite{gennery1999traversability}.
We model the breakdown rate as
\begin{equation}
\phi(\vx) = \frac{\rho(\vx)^2}{5000}
\end{equation}
and it is shown in Figure \ref{fig:marsExample}\subfig{c}.
This scaling 
is largely arbitrary, and can be modified based on the rover's capabilities.

Figures \ref{fig:marsExample}\subfig{d} and \ref{fig:marsExample}\subfig{e} show $u_2$ and $u_1$ and sample optimal trajectories for Example 4.
In mode 1, the rover is able to traverse the crater's rim by avoiding the most steeply sloped portions.
If a breakdown occurs before it has finished climbing, it may become optimal to take a much longer route that avoids the steep slopes of the rim.
However, if a breakdown occurs after the rover has navigated most of the treacherous terrain, it is optimal to climb down the slope towards the target, though with an altered path that avoids steeper slopes.
Figure \ref{fig:marsExample}\subfig{f} shows the elevation of the region of interest, a single mode 1 optimal trajectory, two sample breakdown locations, and subsequent mode 2 paths.

\section{CONCLUSION}
\label{sec:conclusion}

We presented a piecewise-deterministic framework for optimal path planning with two types of random breakdown, total and partial, and showed that each leads to structurally different PDE systems for the value functions.
We also presented an efficient numerical method for solving these systems via value-policy iterations, and verified the first-order accuracy of our approximate solutions experimentally.
Our test problems illustrated the effect that the breakdown type, breakdown rate, and depot/goal placement have on optimal trajectories.
Finally, we computed time-optimal trajectories that account for random partial breakdowns in an environment based on Martian terrain data. 

One obvious extension of our model is to allow the robot to voluntarily terminate mode 2, paying for the full in-place repair before a total breakdown actually occurs. This would in effect replace equation \eqref{eq:u2pde} with a quasi-variational inequality \cite{AndrewsVlad}.
Another natural extension is to consider anisotropic problems, allowing for control-dependent $K$, $f$, $\phi$, and $\lambda$, and taking advantage of the numerical methods in \cite{TsaiChengOsherZhao} or \cite{SethVlad2}.
Extensions could also focus on switching rates, either by allowing them to be time-dependent, or by adding more modes to represent additional levels of damage.
In addition, it will be useful to look beyond the expected cost, maximizing the probability that the cumulative cost falls below any specified threshold \cite{CarteeVlad_UQ}.
Multiobjective extensions are yet another possibility, with hard or chance constraints, similar to \cite{ono2015chance},  imposed on additional criteria (e.g., directly constraining the probability of a full breakdown on Mars).

\addtolength{\textheight}{-0.0cm}

\section*{ACKNOWLEDGMENT}

The authors would like to thank Lars Gr\"{u}ne for his advice on value-policy iteration and Elliot Cartee, whose work on modeling environmental crime inspired this project.


\bibliographystyle{IEEEtran}
\bibliography{IEEEabrv,consolidated_refs}

\end{document}